\documentclass[12pt]{amsart}
\allowdisplaybreaks
\usepackage[english]{babel}
\usepackage[pdftex,paper=a4paper,portrait=true,textwidth=450pt,textheight=675pt,tmargin=3cm,marginratio=1:1]{geometry}
\usepackage{amsfonts}
\usepackage[dvips]{graphics}
\usepackage[colorinlistoftodos]{todonotes}
\usepackage{amsmath}
\usepackage{amsthm}
\usepackage{amssymb}
\usepackage{bbm}
\usepackage{comment}
\usepackage{cancel}
\usepackage{color}
\usepackage[curve]{xypic}
\usepackage{graphicx}
\numberwithin{equation}{section}
\newtheorem{theorem}{Theorem}[section]

\newtheorem{proposition}[theorem]{Proposition}

\newtheorem{conjecture}[theorem]{Conjecture}

\theoremstyle{definition}
\newtheorem{definition}[theorem]{Definition}

\newtheorem{remark}[theorem]{Remark}
\newtheorem{notation}[theorem]{Notation}

\theoremstyle{property}

\DeclareFontFamily{OT1}{rsfs}{}
\DeclareFontShape{OT1}{rsfs}{n}{it}{<-> rsfs10}{}
\DeclareMathAlphabet{\curly}{OT1}{rsfs}{n}{it}

\renewcommand\O{\mathcal O}

\newcommand\EE{\mathbb E}

\newcommand\E{\mathbb E}
\newcommand\C{\mathbb C}

\newcommand\Q{\mathbb Q}

\newcommand\R{\mathbb R}

\newcommand\Z{\mathbb Z}

\newcommand\Coeff{\mathrm{Coeff}}

\newcommand\vd{\mathrm{vd}}
\newcommand\pt{\mathrm{pt}}
\newcommand\vir{\mathrm{vir}}

\newcommand\SW{\mathrm{SW}}

\newcommand\rk{\operatorname{rk}}

\newcommand\ch{\operatorname{ch}}

\newcommand\Pic{\operatorname{Pic}}

\newcommand\INTO{\ar@{^{(}->}[r]}

\newcommand{\smfr}[2]{\hbox{$\frac{#1}{#2}$}}

\DeclareRobustCommand{\SkipTocEntry}[4]{}
\def\eps{\varepsilon}

\begin{document}
\title[Virtual blowup formulas]{Blowup formulas for Segre and Verlinde numbers of surfaces and higher rank Donaldson invariants}
\author[G\"ottsche]{L.~G\"ottsche}
\maketitle

\vspace{-1cm}

\begin{abstract}
We formulate conjectural blowup formulas for Segre and Verlinde numbers on moduli spaces of sheaves on projective surfaces $S$ with 
$p_g(S)>0$ and $b_1(S)=0$. As applications we give a give a conjectural formula for the Donaldson invariants of $S$ in arbitrary rank, 
as well as for the $K$-theoretic Donaldson invariants, and some Donaldson invariants with fundamental matters.

\end{abstract}
\thispagestyle{empty}
\tableofcontents

\section{Introduction} 
\addtocontents{toc}{\protect\setcounter{tocdepth}{1}}
In this paper we  study enumerative invariants of moduli spaces of sheaves on projective algebraic surfaces $S$, with first Betti number $b_1(S)=0$ and geometric genus $p_g(S)>0$. 
Our aim in this paper is to formulate conjectural blowup formulas relating these invariants for the surface $S$ and its blowup $\widehat S$ in a point. 
We will give evidence for these conjectures and then explore some of their consequences. 

Let $H$ be an ample line bundle on $S$.
Let $M:=M_S^H(\rho,c_1,c_2)$ be the moduli space of Gieseker $H$-semistable torsion-free sheaves $E$ on $S$ of rank $\rk(E)=\rho$ and with Chern classes $c_1(E)=c_1$, $c_2(E)=c_2$. 
We denote by $K_0(X)$ the Grothendieck group of coherent sheaves on a scheme $X$, and by $K^0(X)$ the Grothendieck group of locally free sheaves.
$M$ carries a perfect obstruction theory as defined in \cite{Moc} of expected dimension 
$$\vd=\vd(M)=\vd(\rho,c_1,c_2)=2\rho c_2-(\rho-1)c_1^2-(\rho^2-1)\chi(\O_S),$$
and therefore has a virtual fundamental class $[M]^\vir\in H_{2\vd}(M,\Z)$ and a virtual structure sheaf $\O_M^{\vir}$ in $K_0(M)$. To a class $V\in K^0(M)$ we can  associate the virtual holomorphic Euler characteristic
$\chi^\vir(M,V):=\chi(M,V\otimes \O_M^\vir)$. 
We can use this virtual structure to define and study two types of enumerative invariants of $M$, both associated to classes $\alpha\in K^0(S)$:
\begin{enumerate}
\item the virtual Verlinde numbers $\chi^\vir(M,\lambda(\alpha))$, for $\lambda(\alpha)\in \Pic(M)$ the determinant line bundle associated to $\alpha$,
\item the virtual Segre invariants $\int_{[M]^\vir} c(\tau(\alpha))$, for $\tau(\alpha)\in K^0(M)$ a tautological "vector bundle" associated to $\alpha$, which generalizes the tautological bundles $\alpha^{[n]}$ on  Hilbert schemes of points.
\end{enumerate}
We briefly review definitions and properties, both known and conjectured  of these invariants, mostly based on \cite{GKVer}. They generalize the corresponding results on Hilbert schemes of points from \cite{EGL}, \cite{MOP1}, \cite{MOP2}, \cite{MOP3}. More details can be found in Section 2 of this paper and  in \cite{GKVer}.

In this whole paper let $S$ be a smooth projective surface with $b_1(S)=0$, $p_g(S)>0$, and let $H$ be an ample line bundle on $S$. For any $\rho \in \Z_{>0}$, $c_1 \in H^2(S,\Z)$ and $c_2 \in H^4(S,\Z)$, let $M:=M_S^H(\rho,c_1,c_2)$ be the moduli space of rank $\rho$ Gieseker $H$-semistable torsion free sheaves on $S$ with Chern classes $c_1,c_2$. We will always assume in the following that $M$ contains no strictly semistable sheaves. For simplicity of exposition  we also assume that there exists a universal sheaf $\E$ on $S \times M$, but in \cite[Rem.~2.4,~2.5]{GKVer} it is explained how to remove this assumption.  Let $\pi_S : S \times M \rightarrow S$ and $\pi_M : S \times M \rightarrow M$ be the projections to the factors.

{\bf Virtual Segre invariants:}
 Consider the slant product
$$
/ : H^p(S \times M,\Q) \times H_q(S,\Q) \rightarrow H^{p-q}(M,\Q).
$$
 For any  $\gamma \in H^k(S,\Q)$  let
\begin{equation*} 
\mu(\gamma):=  \Big(c_2(\E)-\frac{\rho-1}{2\rho} c_1(\E)^2 \Big) / \mathrm{PD(\gamma )} \in H^{k}(M,\Q),
\end{equation*}
where $\mathrm{PD(\gamma})$ denotes the Poincar\'e dual.
For any class $\alpha \in K^0(S)$, we define 
$$
\ch(\tau(\alpha)) :=  \ch(- \pi_{M !} ( \pi_S^* \alpha \cdot \E \cdot  \det(\E)^{-\frac{1}{\rho}})) \in A^*(M)_{\mathbb{Q}},
$$
where $A^*(M)_{\Q}$ denotes the Chow ring with rational coefficients. 
Let $\pt \in H^4(S,\Z)$ be the Poincar\'e dual of the point class, and $L\in H^2(S,\Z)$.
The virtual Segre invariant of $S$ associated to $(\rho,c_1,c_2,\alpha,L)$ is 
$\int_{[M]^{\vir}} c(\tau( \alpha)) \, \exp(\mu(L)  + \mu(\pt) u)$, where $u$ is a formal variable.
In the future we will always put $s:=\rk(\alpha)$.

{\bf Virtual Verlinde invariants:} Let $c \in K(S)_{\mathrm{num}}$ be a class in the numerical Grothendieck group of $S$ such that $\rk(c) = \rho$, $c_1(c) = c_1$, and $c_2(c) = c_2$. Denote
\begin{equation*} 
K_c := \{ v \in K^0(S) \, : \, \chi(S,c \otimes v) = 0\}.
\end{equation*} 
The map
\begin{equation} \label{lambdaE}
\lambda:=\lambda_{\EE} : K_c \rightarrow \Pic(M), \quad \beta \mapsto \det \big( \pi_{M!} \big( \pi_S^* \beta \cdot [\E] \big) \big)^{-1}.
\end{equation}
does not depend on the choice of the universal family and can be extended to the case where the universal family does not exist.
Fix $r \in \Z$, $L \in \Pic(S) \otimes \Q$ such that $\mathcal{L}:=L \otimes \det(c)^{-\frac{r}{\rho}} \in \Pic(S)$ and  $\rho$ divides $\mathcal{L}c_1 + r \big( \frac{1}{2} c_1(c_1-K_S) - c_2 \big)$. Take $v \in K^0(S)$ such that
\begin{itemize}
\item $\rk(v) = r$ and $c_1(v) = \mathcal{L}$, 
\item $c_2(v) = \frac{1}{2} \mathcal{L}(\mathcal{L}-K_S) +r\chi(\O_S) + \frac{1}{\rho} \mathcal{L} c_1 + \frac{r}{\rho}\big(\frac{1}{2}c_1(c_1-K_S) - c_2\big)$.
\end{itemize}
The second condition is equivalent to $v \in K_c \subset K^0(S)$. We define
$
\mu(L) \otimes E^{\otimes r} := \lambda(v).
$
Let $\O_M^\vir$ be the virtual structure sheaf of $M$. For $V\in K^0(M)$ let $\chi^\vir(M,V):=\chi(M,V\otimes \O_M^\vir)$.
The virtual Verlinde number of $S$ associated to $(\rho,c_1,c_2,L,r)$ is 
$\chi^\vir(M,\mu(L) \otimes E^{\otimes r})$.

If $\rho=1$ and $c_1=0$, $M_S^H(1,0,n)$ is just the Hilbert schemes $S^{[n]}$ of $n$ points on $S$. 
Let $Z_n(S)\subset S\times S^{[n]}$ be the universal subscheme with the projections $p:Z_n(S)\to S^{[n]}$ and $q:Z_n(S)\to S$.
For $\alpha\in K^0(S)$, the associated tautological bundle on $S^{[n]}$ is defined by $\alpha^{[n]}=p_*q^*(\alpha)\in K^0(S^{[n]}$. Then we have
$\alpha^{[n]}=\tau(\alpha)$ in the notations above, thus the virtual Segre invariants generalize the Segre numbers
$\int_{S^{[n]}} c(\alpha^{[n]})$ on Hilbert schemes of points. Similarly if $r=\rk(\alpha)$ and $L=\det(\alpha)$, we have that 
 $\mu(L)\otimes E^{r})=\det(\alpha^{[n]})$, thus the virtual Verlinde invariants generalize the Verlinde numbers of line bundles on Hilbert schemes of points.
 
\subsection*{Conjectural formulas and relations}
We state the conjectures of \cite{GKVer} expressing  the Segre and Verlinde numbers in terms of Seiberg-Witten invariants.
These are conjectural generalizations of results in  \cite{EGL}, \cite{MOP1}, \cite{MOP2}, \cite{MOP3} for Hilbert schemes of points.
We denote the Seiberg-Witten invariant of $S$ in class $a \in H^2(S,\Z)$
by $SW(a).$ We follow Mochizuki's  \cite{Moc} convention $SW(a) = \widetilde{SW}(2a-K_S)$, where $\widetilde{SW}$ denotes the usual Seiberg-Witten invariant  from differential geometry. For $S$ a surface with $b_1(S)=0$ and $p_g(S)>0$,  there are finitely many $a \in H^2(S,Z)$ with  $SW(a)\ne 0$, called Seiberg-Witten basic classes. 

\begin{notation}
We fix $\rho\in \Z_{>0}$, and $s,r\in \Z$.
We denote $\eps_\rho := \exp(2 \pi i /\rho)$ and write $[n]:=\{1, \ldots, n\}$ for any $n \in \mathbb{Z}_{\geq 0}$.  For integers $n<m$ we also write $[n,m]$ for the integers $l$ 
with $n\le l\le m$. 
We introduce the following changes of variables 
$$z = t (1+(1-\tfrac{s}{\rho}) t^2)^{\frac{1}{2}(1-\frac{s}{\rho})},\quad w=v(1+v^2)^{\frac{1}{2}(\frac{r^2}{\rho^2}-1)},$$ 
and the following formal power series
\begin{align*}
V^{(\rho)}_s(z) &= (1+(1-\tfrac{s}{\rho})t^2)^{\rho-s} (1+(2-\tfrac{s}{\rho})t^2)^s, \\ 
W^{(\rho)}_s(z) &= (1+(1-\tfrac{s}{\rho})t^2)^{\frac{1}{2}(s-\rho-1)} (1+(2-\tfrac{s}{\rho})t^2)^{\frac{1}{2}(1-s)}, \\
X^{(\rho)}_s(z) &=  (1+(1-\tfrac{s}{\rho})t^2)^{\frac{1}{2}(s^2-(\rho+\frac{1}{\rho})s)} (1+(2-\tfrac{s}{\rho})t^2)^{-\frac{1}{2}s^2+\frac{1}{2}} (1+(1-\tfrac{s}{\rho})(2-\tfrac{s}{\rho})t^2)^{-\frac{1}{2}}, \\
Q^{(\rho)}_s(z) &= \tfrac{1}{2}t^2(1+(1-\tfrac{s}{\rho}) t^2), \quad R^{(\rho)}_s(z) = t^2, \quad T^{(\rho)}_s(z) = \rho t^2 (1+ \tfrac{1}{2}(1-\tfrac{s}{\rho})(2-\tfrac{s}{\rho})t^2),\\
G_r(w)&= 1+v^2, \quad 
F^{(\rho)}_r(w) = \frac{(1+v^2)^{\frac{r^2}{\rho^2}}}{(1+\tfrac{r^2}{\rho^2} v^2)}.
\end{align*}
\end{notation}

We warn the reader that this notation differs from that of  \cite{GKVer}: $z^{\frac{1}{2}},t^{\frac{1}{2}},w^{\frac{1}{2}},v^{\frac{1}{2}}$ from \cite{GKVer} were replaced by where $z,t,w,v$. 
Also in the following conjectures $Z^{(\rho)}_s$ from  \cite{GKVer} was replaced by $\rho Z^{(\rho)}_s$ and $B^{(\rho)}_r$ by $\rho B^{(\rho)}_r$,

\begin{conjecture} \label{conj2}\cite[Conj.~2.8,~2.9]{GKVer}
For $\rho  \in \Z_{>0}$ and $s,r \in \Z$, there exist $Y^{(\rho)}_s$, $Z^{(\rho)}_s$, $Y^{(\rho)}_{j,s}$, $Z^{(\rho)}_{jk,s}$, $S^{(\rho)}_s$, $S^{(\rho)}_{j,s} \in \C[[z]]$ and $A^{(\rho)}_{r}$, $B^{(\rho)}_{r}$, $A^{(\rho)}_{j,r}$, $B^{(\rho)}_{jk,r} \in \C[[w]]$ for all $1 \leq j \leq k \leq \rho-1$ such that the following holds. For any $\alpha \in K^0(S)$ with $\rk(\alpha)=s$ and $L \in \Pic(S)$, we put 
\begin{align*}
\phi_{S,\rho,c_1,\alpha,L}&:=\rho^{2 - \chi(\O_S)} \, V_s^{c_2(\alpha)} \, W_s^{c_1(\alpha)^2} \, X_s^{\chi(\O_S)} \, Y_{s}^{c_1(\alpha) K_S} \, Z_{s}^{K_S^2} \, e^{L^2 Q_s + (c_1(\alpha)L) R_s + (LK_S) S_{s}  + u \, T_s} \\
&\quad \cdot \sum_{(a_1, \ldots, a_{\rho-1})} \prod_{j=1}^{\rho-1} \eps_{\rho}^{j a_j c_1} \, \SW(a_j) \, Y_{j,s}^{c_1(\alpha) a_j} \, e^{(a_j L) S_{j,s}} \prod_{1 \leq j \leq k \leq \rho-1} Z_{jk,s}^{a_j a_k},\\
\psi_{S,\rho,c_1,L,r}&:=\rho^{2 - \chi(\O_S)} \, G_{r}^{\chi(L)} \, F_{r}^{\frac{1}{2} \chi(\O_S)} \, A_{r}^{LK_S} \, B_{r}^{K_S^2} \\ &\quad \cdot \sum_{(a_1, \ldots, a_{\rho-1})} \prod_{j=1}^{\rho-1} \eps_{\rho}^{j a_j c_1} \, \SW(a_i) \, A_{j,r}^{a_j L} \prod_{1 \leq j \leq k \leq \rho-1} B_{jk,r}^{a_j a_k}.
\end{align*}
Let  $M=M_S^H(\rho,c_1,c_2)$ consist only of stable sheaves. Then
\begin{align*}
\int_{[M]^{\vir}} c( \tau(\alpha) ) \, \exp(\mu(L)  + \mu(\pt) u)=\Coeff_{z^{\vd(M)}}[\phi_{S,\rho,c_1,\alpha,L}],\\
\chi^{\vir}(M, \mu(L) \otimes E^{\otimes r})=\Coeff_{w^{\vd(M)}}[\psi_{S,\rho,c_1,L,r}].
\end{align*}
where the sums are over all $(a_1, \ldots, a_{\rho-1}) \in H^2(S,\Z)^{\rho-1}$. 
Furthermore, $Y_s$, $Y_{j,s}$,  $Z_s$, $Z_{jk,s}$, $S_s$, $S_{j,s}$, $A^{(\rho)}_{r}$, $B^{(\rho)}_{r}$, $A^{(\rho)}_{j,r}$, $B^{(\rho)}_{jk,r}$ are all algebraic functions.
\end{conjecture}

In this paper  we  formulate conjectural blowup formulas for the Verlinde and Segre invariants of moduli spaces of sheaves on surfaces. These formulas give a simple relationship between the invariants on a surface $S$ and some of the invariants on its blowup $\widehat S$ in a point. 
For moduli spaces $M$ of sheaves of rank $\rho$, most of them hold for the Verlinde numbers $\chi^\vir(M,\mu(L)\otimes E^{\otimes r})$ under the condition that 
$|r|\le \rho$,   and for the Segre invariants  under the condition that $0\le s\le 2\rho$.
We have found these relations experimentally as relations among the universal functions $A_{J,r}^{(\rho)}$, $B_{J,r}^{(\rho)}$ for the Verlinde numbers and 
$Y_{J,r}^{(\rho)}$, $Z_{J,r}^{(\rho)}$, $S_{J,r}^{(\rho)}$ for the Segre invariants, which we will introduce in Section 2. Here we will state them directly for the generating functions of the invariants of the moduli spaces, where they
take a simpler and more attractive form. We start with the blowup formulas for the Verlinde numbers.

\begin{conjecture}\label{blowpsi}
For any smooth projective surface $S$ with 
$b_1(S)=0$ and $p_g(S)>0$,
 let $\pi:\widehat S \to S$ be the blowup in a point with exceptional divisor $D$. Let $|r|\le \rho$.
\begin{enumerate}
\item 
For $k\in[0,\rho]$ and in addition, if $|r|<\rho$, for $k\in [r,r+\rho]$ we have 
\begin{align}\label{psiblow1}
\psi_{\widehat S,\rho,\pi^*c_1,L+kD,r}&=\psi_{S,\rho,c_1,L,r}, \quad k\in [0,\rho].
\end{align}
\item For $\ell=1,\ldots,\rho-1$ we have 
\begin{equation}\label{psiblow2}
\psi_{\widehat S,\rho,\pi^*c_1-\ell D,\pi^*L+(k+(\rho-\ell)\frac{r}{\rho})D,r}
=\begin{cases}0,&k\in[1,\rho-1],\\
\rho\cdot (-w)^{\ell(\rho-\ell)}\psi_{S,\rho,c_1,L,r},& k=0,\  r\in[-\rho+1,\rho],\\
\rho\cdot (-w)^{\ell(\rho-\ell)}\frac{\psi_{S,\rho,c_1,L,r}}{(1-(-v^2)^\rho)^{\rho-\ell}},& k=0,\ r=-\rho,\\
\rho\cdot w^{\ell(\rho-\ell)}\psi_{S,\rho,c_1,L,r},&k=\rho,\ r\in [-\rho,\rho-1],\\
\rho\cdot w^{\ell(\rho-\ell)}\frac{\psi_{S,\rho,c_1,L,r}}{(1-(-v^2)^\rho)^{\rho-\ell}},& k=\rho,\ r=\rho.
\end{cases}
\end{equation}
\end{enumerate}
\end{conjecture}

For the Segre invariants, we get the following blowup formulas.

\begin{conjecture} \label{blowphi}
Let $S$ be a surface with $p_g(S)>0$, $b_1(S)=0$. Let $L\in \Pic(S)$, let $\widehat S$ be the blowup of $S$ in a point with exceptional divisor $D$. Let $\alpha\in K^0(S)$ be a class of rank $s\in [0,2\rho].$ Let $L\in H^2(S,\Q)$.
\begin{enumerate}
\item For $k=[0,\rho]$, and furthermore, if $s\in [1,2\rho-1]$ for $k\in [s-\rho,s]$, we have
$$\phi_{\widehat S,\rho,\pi^*c_1,\pi^*\alpha+kO_D(D),\pi^*L}=\phi_{S,\rho,c_1,\alpha,L}.$$
\item For $k\ge 0$,
$$\phi_{\widehat S,\rho,\pi^*c_1,\pi^*\alpha+kO_D(D),\pi^*L-xD}= \phi_{S,\rho,c_1,\alpha,L}+ O(x^{\rho-k+1}),$$
 and furthermore, if $s\in [1,2\rho-2]$, we have for $k\le s$
 $$\phi_{\widehat S,\rho,\pi^*c_1,\pi^*\alpha+kO_D(D),\pi^*L-xD}= \phi_{S,\rho,c_1,\alpha,L}+ O(x^{\rho+k-s+1}),$$
\item Let $\widehat c_1\in H^2(\widehat S,\Z)$ with $\widehat c_1D=\ell\in [1,\rho-1]$, let $\widehat \alpha\in K^0(\widehat S)\otimes \Q$ with $c_1(\widehat \alpha)D=\ell\frac{s}{\rho}-s-k$ with $k\in \Z$. Then 
$$\phi_{\widehat S,\rho,\widehat c_1 ,\widehat \alpha,\pi^*L-xD}=O(x^{\min(l(\rho-l+k),(l-k)(\rho-l))}).$$
In particular $\phi_{\widehat S,\rho,\widehat c_1 ,\widehat \alpha,\pi^*L}=0$ for $k\in [-\rho+\ell+1,\ell-1]$.
\end{enumerate}
\end{conjecture}

Together with the known and conjectured properties of the Verlinde and Segre invariants mentioned above,  these conjectural  blowup formulas impose rather strong constraints. In many cases this  allows  to determine conjectural formulas for the universal functions $A^{(\rho)}_{J,r}$, $B^{(\rho)}_{J,r}$, $Y^{(\rho)}_{J,s}$, $Z^{(\rho)}_{J,s}$, $S^{(\rho)}_{J,s}$. In \cite{Gstr} we first carry this out  in case  $|r|\le \rho$  or $0\le s\le 2\rho$ for $\rho\le 6$,  then using also a virtual version of strange duality give a conjectural formula for  the $A^{(\rho)}_{J,r}$ and  $Y^{(\rho)}_{J,s}$ for arbitrary values of $\rho$ and $r$ (respectively $s$).

{\bf Donaldson invariants in arbitrary rank.} In this paper we will concentrate on the  values $r=-\rho,0,\rho$ for the Verlinde numbers  and $s=0,\rho,2\rho$ for the Segre invariants.
In these cases we will be able to determine  conjectural generating functions for the Segre and Verlinde invariants of the moduli spaces $M=M_S^H(\rho,c_1,c_2)$ of sheaves of arbitrary rank $\rho$. Of particular interest are the Donaldson invariants of arbitrary rank $\rho$, which are a special instance of the case $s=0$ of the Segre invariants.
The Segre invariants with a given $s$ are sometimes in the physics language called Donaldson invariants with $s$ fundamental matters. Thus we conjecturally determine
also the Donaldson invariants in rank $\rho$ with $\rho$ or $2\rho$ fundamental matters. The Verlinde numbers with $r=0$ are sometimes called $K$-theoretic Donaldson invariants.

Donaldson invariants are invariants of differentiable $4$-manifolds $X$ defined using moduli spaces of anti-selfdual connections on principal $SU(2)$ and $SO(3)$ bundles on $X$. If $X$ is a projective algebraic surface $S$,  they can be computed as intersection numbers of $\mu$-classes on moduli spaces of sheaves $M_S^H(2,c_1,c_2)$ on $S$.  The definition of Donaldson invariants has been generalized in \cite{Kro} to higher rank bundles. Again they can be computed  using moduli spaces of higher rank sheaves $M_S^H(\rho,c_1,c_2)$.

\begin{definition} 
Let $(S,H)$ be a polarized surface and fix $\rho\in \Z_{>0}$, $c_1\in H^2(S,\Z)$. Let $L\in H^2(S,\Q)$. The rank $\rho$ Donaldson invariants of $S$ with respect to $H,c_1$ are 
$$D^{S,H}_{\rho,c_1,c_2}(L  +u \pt)=\int_{[M_S^H(\rho,c_1,c_2)]^\vir} \exp(\mu(L) + \mu(\pt)u).$$
\end{definition}

For $a\in H^2(S,\Z)$ we put $\widetilde a:=2a-K_S\in H^2(S,\Z)$. As mentioned above, we  denote by $\widetilde{SW}$ the Seiberg-Witten invariants in gauge theory, which are related to our convention by $\widetilde{SW}(\widetilde a)=SW(a)$. The Witten conjecture \cite{Wit} expresses the  rank $2$ Donaldson invariants in terms of the Seiberg-Witten invariants. It was proven in \cite{GNY3} for projective surfaces, and (modulo a technical assumption) for all differentiable $4$-manifolds $M$ with $b_1(M)=1$, odd $b_+(M)$ and of Seiberg-Witten simple type in \cite{FL1,FL2}.
The Marino-Moore conjecture \cite{MM} extends the Witten conjecture to Donaldson invariants of arbitrary rank.  The following conjecture can be viewed as a completely explicit version of the Marino-Moore conjecture.

\begin{notation} In the whole paper we let $\xi=e^{\frac{\pi i}{2\rho}}$ be a primitive $4\rho$-th root of unity.  For integers $n,m$ we denote $[n,m]:=\bigr\{k\in \Z\bigm| n\le k\le m\bigr\}$, and we put $[m]:=[1,m]$.
For $1\le i\le j\le\rho-1$ let 
$\beta_{ij}:=\frac{\xi^{i+j}-\xi^{-(i+j)}}{\xi^{j-i}-\xi^{i-j}}$ and $\beta_{ji}:=\beta_{ij}$.  For any subset $J\subset [\rho-1]$ put 
$$\beta_J:=\prod_{\substack{i\in J\\j\in [\rho-1]\setminus J}}\beta_{ij},$$
and let 
$$B:=\sum_{J\subset [\rho-1]} \beta_{J}.$$
\end{notation}

\begin{conjecture} \label{Donconj}
Let $(S,H)$ be a polarized surface and fix $\rho\in \Z_{>0}$, $c_1\in H^2(S,\Z)$, $c_2\in \Z$, and write $\vd:=\vd(\rho,c_1,c_2)$. We also write $\xi=e^{\frac{\pi i}{2\rho}}$, and $\varepsilon_\rho=e^{\frac{2\pi i}{\rho}}$. Then 
$D^{S,H}_{\rho,c_1,c_2}(L  +u \pt)$ is the coefficient of $z^{\vd}$ of
\begin{align}\label{Donform}
\rho^{2-\chi(\O_S)}B^{K_S^2}e^{(\frac{1}{2}L^2+\rho u)z^2}\sum_{(a_1,\ldots,a_{\rho-1})}\prod_{j=1}^{\rho-1}\varepsilon_\rho^{j\cdot(a_j,c_1)}\widetilde{SW}(\widetilde a_j)
e^{-\sin(\pi\frac{j}{\rho})(\widetilde a_j L)z}\prod_{1\le i< j\le \rho-1}\beta_{ij}^{\frac{1}{2}\widetilde a_i(\widetilde a_j-\widetilde a_i)}
\end{align}
Note that by definition for $1\le i<j\le \rho-1$, we can write $\beta_{ij}=\frac{\sin(\pi\frac{i+j}{\rho})}{\sin(\pi \frac{j-i}{\rho})}\in \R_{>0}$, and we choose $\beta_{ij}^{\frac{1}{2}}\in \R_{>0}$.
\end{conjecture}
This extends  \cite[Conj.~5.6]{GKVer}, where the shape of the formula for general $\rho$ is predicted and the explicit formula given for $\rho\le 4$. 

Let $M=M_S^{H}(\rho,c_1,c_2)$, and assume that $\vd(\rho,c_1,c_2)=0$. Then Conjecture \ref{Donconj} in particular gives a  formula for 
the degree of $[M]^\vir$: 
$$\int_{[M]^\vir} 1=\rho^{2-\chi(\O_S)}B^{K_S^2}\sum_{(a_1,\ldots,a_{\rho-1})}\prod_{j=1}^{\rho-1}\varepsilon_\rho^{j\cdot(a_j,c_1)}\widetilde{SW}(\widetilde a_j)
\prod_{1\le i< j\le \rho-1}\beta_{ij}^{\frac{1}{2}\widetilde a_i(\widetilde a_j-\widetilde a_i)}.$$
In \cite{GKL} we compute generating functions of the monopole contributions to Vafa-Witten invariants using localization on Hilbert schemes of points for $\rho\le 5$, then using the $S$-duality
conjecture we get conjectural  generating functions for the virtual Euler numbers of the moduli spaces $M_S^H(\rho,c_1,c_2)$. In particular this gives us conjectural formulas for the degrees of the virtual fundamental classes of moduli spaces $M_S^H(\rho,c_1,c_2)$ of virtual dimension zero for $\rho\le 5$. These formulas also agree with the formulas of Conjecture \ref{Donconj}.

\noindent \textbf{Acknowledgements.} This work grew out of my collaboration with Martijn Kool over the last few years, in which I benefitted from 
numerous discussions and insights, without which the current paper could not have been written. The results and ideas of  \cite{GKVer} play an important role in this paper.

\section{Background material}
In this section we briefly review some of the conjectures and definitions from \cite{GKVer}, which can be viewed 
as generalizations of corresponding conjectures and results from \cite{EGL}, \cite{MOP1}, \cite{MOP2}, \cite{MOP3} from  the case of Hilbert schemes of points to moduli spaces of sheaves of higher rank. For more details see \cite{GKVer}. 

\begin{notation}
It is sometimes convenient to organize the universal functions of Conjecture \ref{conj2} in a different form.
We put  
\begin{align}
\begin{split} \label{Jseries}
Y_{J,s}&:=Y_s\prod_{j\in J}Y_{j,s}, \quad Z_{J,s}:=Z_s\prod_{i\le j\in J}Z_{ij,s},\quad S_{J,s}:=S_s+\sum_{j\in J}S_{j,s}, \\
 A_{J,r}&:=A_r\prod_{j\in J}A_{j,r}, \quad B_{J,r}:=B_r\prod_{i\le j\in J}B_{ij,r}.
 \end{split}
 \end{align}
 We note that knowing the power series $Y_{J,s}$, $Z_{J,s}$, $S_{J,s}$, $A_{J,r}$, $B_{J,r}$ is equivalent to knowing the power series
$Y_s$, $Y_{i,s}$, $Z_s$, $Z_{ij,s}$, $S_s$, $S_{i,s}$, $A_r$, $A_{i,r}$, $B_r$, $B_{ij,r}$. 

For a subset $|J|\subset [\rho-1]$ we also let $|J|$ be the number of elements of $J$ and put $\|J\|:=\sum_{j\in J} j$.
\end{notation}
Let $S$ have an irreducible  canonical divisor. Then  it is well known that the Seiberg-Witten classes of $S$  are just $0$ and $K_S$, with Seiberg-Witten invariants 
$SW(0)=1$, $SW(K_S)=(-1)^{\chi(\O_S)}$. Thus Conjecture \ref{conj2} gives the following.
\begin{conjecture}\label{conj1}\cite[Conj.~1.4,~1.6]{GKVer}
\begin{enumerate}
\item 
For any $\alpha \in K^0(S)$ with $\rk(\alpha)=s$ and any $L \in \Pic(S)$,   we have that  $\int_{[M]^{\vir}} c( \tau(\alpha)) \, \exp(\mu(L)  + \mu(\pt) u)$ is the coefficient of $z^{\vd(M)}$ of
\begin{align*}
\rho^{2 - \chi(\O_S)} \, V_s^{c_2(\alpha)}& \, W_s^{c_1(\alpha)^2} \, X_s^{\chi(\O_S)} \, e^{L^2 Q_s + (c_1(\alpha)L) R_s + u \, T_s} \\&\quad \cdot\sum_{J \subset [\rho-1]} (-1)^{|J| \chi(\O_S)}  \, \eps_{\rho}^{\|J\| K_S c_1} \, Y_{J,s}^{c_1(\alpha) K_S} \, Z_{J,s}^{K_S^2} \, e^{(K_S L) S_{J,s}}.
\end{align*}
\item 
The virtual Verlinde number $\chi^{\vir}(M, \mu(L) \otimes E^{\otimes r})$ equals the coefficient of $w^{\vd(M)}$ of
\begin{align*}
\rho^{2 - \chi(\O_S)} \, G_{r}^{\chi(L)} \, F_{r}^{\frac{1}{2} \chi(\O_S)} \sum_{J \subset [\rho-1]} (-1)^{|J| \chi(\O_S)} \, \eps_{\rho}^{\|J\| K_S c_1} \, A_{J,r}^{K_S L} \, B_{J,r}^{K_S^2}.
\end{align*}
\end{enumerate}
Furthermore, the $Y_{J,s}$,   $Z_{J,s}$,  $S_{J,s}$,  $A_{J,r}$, $B_{J,r}$ are all algebraic functions.
\end{conjecture}

The Segre-Verlinde correspondence, conjectured first in the case of Hilbert schemes of points in \cite{Joh} and \cite{MOP2}
gives an explicit conjectural relation between the universal functions of the Verlinde and the Segre numbers. 
\begin{conjecture} \label{SegreVerlindeconj} \cite[Conj.~1.7]{GKVer} For any $\rho>0$, $r \in \Z$, $s:=\rho+r$ we have 
\begin{align*}
A_{J,r}(w) = W_{s}(z) \, Y_{J,s}(z), \quad  B_{J,r}(w) = Z_{J,s}(z), 
\end{align*}
for all $J \subset [\rho-1]$, under the change of variables
\begin{equation} \label{varchange2}
w = v(1+v^2)^{\frac{1}{2}(\frac{r^2}{\rho^2}-1)}, \quad z = t (1+(1-\tfrac{s}{\rho}) t^2)^{\frac{1}{2}(1-\frac{s}{\rho})}, \quad v = t( 1- \tfrac{r}{\rho} t^2 )^{-\frac{1}{2}}.
\end{equation}
\end{conjecture}

Finally the virtual Serre duality $\chi^\vir(M,L)=(-1)^{\vd(M)}\chi^\vir(M,L^\vee\otimes K_M^\vir)$ (see \cite{FG}), leads to the following conjecture.
\begin{conjecture} \label{SerreDuality}\cite[Conj.~5.4]{GKVer}
For any $\rho > 0$, we have
\begin{align*}
A_{J,-r}(w)&=(1+v^2)^{1-\rho} A_{J,r}(-w)^{-1}, \\
B_{J,-r}(w)&=(1+v^2)^{\binom{\rho}{2}}A_{J,r}(-w)^\rho B_{J,r}(-w),
\end{align*}
for all $J \subset [\rho-1]$ and $r \in \Z$.
\end{conjecture}

In \cite[Conj.~3.5]{GKVer} it is also conjectured that the coefficients of the universal functions $A_{J,r}, B_{J,r}$ $Y_{J,s}$, $Z_{J,s}$,  $S_{J,s}$ are polynomials in $r$ and $s$ respectively. We can slightly strengthen the statement there.

\begin{conjecture} \label{poly} For all $J\in [\rho-1]$ and all $n\ge 0$ we have the following. 
\begin{enumerate} \item 
 The coefficients of $w^n$ of  $A_{J,r},$ (respectively $B_{J,r}$) are polynomials in $r$ of degree at most $n-1$ (respectively at most $n$).
\item The coefficients of $z^n$ of  $Y_{J,s},$ (respectively $Z_{J,s},$  $S_{J,s}$) are polynomials in $s$ of degree at most $n-1$ (respectively at most $n$). 
\end{enumerate}
\end{conjecture}

We recall the blowup relation for the virtual Verlinde numbers  \cite[Prop.~5.10]{GKVer} and state an analogous result for the Segre invariants, 
whose proof is an easy adaptation of the proof of the formula for Verlinde numbers.
For this let $S$ be a smooth projective surface satisfying $b_1(S)=0$ and $p_g(S)>0$. 
Let $\pi:\widehat S\to S$ be the blowup of $S$ in a point with exceptional divisor $D$.

\begin{proposition}\label{verblowprop}
\begin{align}\label{verblowform}
\psi_{\widehat{S},\rho,\pi^*c_1-\ell D,L-mD,r} &= \frac{\psi_{S,\rho,c_1,L,r}}{(1+v^2)^{\binom{m+1}{2}}}  \sum_{J \subset [\rho-1]}  \eps_{\rho}^{\|J\| l} \frac{A_{J,r}^{m}}{B_{J,r}}  ,\\
\label{segreblowform}
 \begin{split}
\phi_{\widehat{S},\rho,\pi^* c_1-\ell D,\alpha+k\O_D(D),L-xD} &=\frac{e^{-x^2Q_s+kxR_s}(1+(1-\frac{s}{\rho})t^2)^{k(k-s+\rho)/2}\phi_{S,\rho,c_1,\alpha,L}  }{(1+(2-\frac{s}{\rho})t^2)^{k(k-s)/2}}\\&\qquad\qquad \cdot \sum_{J \subset [\rho-1]}  \eps_{\rho}^{\|J\| \ell} \frac{Y_{J,s}^{k}}{{Z_{J,s}}} e^{xS_{J,s}}.
\end{split}
\end{align}
\end{proposition}

To get the explicit formula \eqref{segreblowform}, note that 
\begin{align*}
c_1(\alpha+k\O_D(D))^2&=c_1(\alpha)-k^2, \quad c_2(\alpha+k\O_D(D))=c_2(\alpha)-\binom{k}{2},\\
V_s^{\binom{k}{2}} W_s^{k^2}&= (1+(1-\smfr{s}{\rho})t^2)^{-\frac{1}{2}k(k-s+\rho)}(1+(2-\smfr{s}{\rho})t^2)^{\frac{1}{2}k(k-s)}.
\end{align*}

\section{Virtual blowup formulas}
We now state the blowup formulas that we experimentally found for the universal functions $A_{J,r}$, $B_{J,r}$ for the Verlinde numbers and 
$Y_{J,s}$, $Z_{J,s}$, $S_{J,s}$ for the Segre invariants. We then show that they are equivalent to  Conjectures \ref{blowphi} and \ref{blowpsi} of the introduction. 

\subsection{Blowup formula for the Verlinde numbers}
\begin{conjecture} \label{blowver1} Let $|r|\le\rho$. With the change of variables $w=v(1+v^2)^{\frac{1}{2}(\frac{r^2}{\rho^2}-1)}$, we have the following.
\begin{enumerate}
\item For $a\in [-\rho,0]$ and in addition, if $|r|<\rho$ for $a\in[-\rho-r,-r] $
\begin{equation}
\label{Ablow1}
\sum_{J\subset [\rho-1]} A_{J,r}^aB_{J,r}^{-1}=(1+v^2)^{\binom{a+1}{2}}.
\end{equation}
\item 
For $\ell=1,\ldots,\rho-1$ , and $a= i+(\ell-\rho)\frac{r}{\rho},$ with  $i\in [-\rho,0]$, we have the following
\begin{equation}\label{Ablow2}
\sum_{J\subset [\rho-1]} \epsilon_\rho^{\ell\|J\|}A_{J,r}^{a}B_{J,r}^{-1}
=\begin{cases}0,&i\in [-\rho+1,-1],\\
(-w)^{\ell(\rho-\ell)}(1+v^2)^{\binom{a+1}{2}},& i=0,\  r\in[-\rho+1,\rho],\\
(-w)^{\ell(\rho-\ell)}\frac{(1+v^2)^{\binom{a+1}{2}}}{(1-(-v^2)^\rho)^{\rho-\ell}},&  i=0,\ r=-\rho,\\
w^{\ell(\rho-\ell)}(1+v^2)^{\binom{a+1}{2}},&i=-\rho,\ r\in [-\rho,\rho-1],\\
w^{\ell(\rho-\ell)}\frac{(1+v^2)^{\binom{a+1}{2}}}{(1-(-v^2)^\rho)^{\rho-\ell}},& i=-\rho,\ r=\rho.
\end{cases}
\end{equation}
\end{enumerate}
\end{conjecture}

We now want to see that, assuming Conjecture \ref{conj2}, the two Conjectures  \ref{blowver1} and  \ref{blowpsi} are equivalent. 
By  Conjecture \ref{conj2}, the universal power series $A_{J,r}$, $B_{J,r}$ are independent of the surface $S$. 
We can therefore assume that $S$ is a K3-surface, and in particular that $\psi_{S,\rho, c_1,L,r}\ne 0$. We put $a=-k$, then we can rewrite 
 \eqref{verblowform} as 
\begin{equation}\label{blowform1}\frac{\psi_{\widehat{S},\rho,\pi^*c_1-\ell D,L+kD,r}}{\psi_{S,\rho,c_1,L,r}} = \frac{1}{(1+v^2)^{\binom{a+1}{2}}} \sum_{J \subset [\rho-1]}  \eps_{\rho}^{\|J\| \ell} \frac{A_{J,r}^{a}}{B_{J,r}}.\end{equation}
Putting $\ell=0$ in \eqref{blowform1}  shows immediately that \eqref{Ablow1} is equivalent to \eqref{psiblow1}, and similarly in case $\ell\ne 0$
\eqref{Ablow2} is equivalent to \eqref{psiblow2}.

We can also directly formulate these conjectures as conjectures about the virtual holomorphic Euler characteristics of moduli spaces.
If $c\in K^0(S)$ is the class of a sheaf $E$ in $M_S^H(\rho,c_1,c_2)$, and  $H$ is an ample line bundle such that $M_S^H(\rho,c_1,c_2)$ only consists of stable sheaves,  we also write $M_S(c):=M_S^H(\rho,c_1,c_2)$.

\begin{conjecture} \label{blowver}
Let $S$ be a surface with $p_g(S)>0$, $b_1(S)=0$. Let $L\in \Pic(S)$, let $\pi:\widehat S\to S$ be the blowup of $S$ in a point with exceptional divisor $D$. Let 
$c\in K^0(S)$ be a class of rank $\rho$, and let $\beta\in K^0(S)$ with  $r:=rk(\beta)$ with $|r|\le \rho$.
\begin{enumerate}
\item For $k\in [0, \rho]$ we have the following:
$$\chi^\vir\big(M_{\widehat S}(\pi^*c),\lambda(\pi^*(\beta)+k\O_D(D))\big)=\chi^\vir\big(M_{S}(c),\lambda(\beta)\big),$$
and if $|r|<\rho$, then also $$\chi^\vir\big(M_{\widehat S}(\pi^*c),\lambda(\pi^*(\beta)+(r+k)\O_D(D))\big)=\chi^\vir\big(M_{S}(c),\lambda(\beta)\big),$$
\item   For $ \ell=1,\ldots,\rho-1$, and $k=1,\ldots,\rho-1$ we have 
$$\chi^\vir\big(M_{\widehat S}(\pi^*c-\ell\O_D),\mu\big(L+(k+\smfr{(\rho-\ell)r}{\rho})D\big)\otimes E^{\otimes r}\big)=0.$$
\end{enumerate}
\end{conjecture}

We see that Conjecture  \ref{blowpsi} implies Conjecture \ref{blowver}: By definition we have
$$\chi^\vir(M_S(c),\mu(L)\otimes E^{\otimes r})=\Coeff_{w^{\vd(M_S(c))}}\big[\psi_{S,\rho,c_1,L,r}\big].$$
To see part (1) of Conjecture \ref{blowver}, note that $\vd(M_S(c))=\vd(M_{\widehat S}(\pi^*c))$, and if $\mu(L)\otimes E^{\otimes r}=\lambda(\beta)$ for some 
$\beta\in K_c(S)$ or rank $r$, then $\pi^*\beta+k\O_D(D)\in K_{\pi^*c}(\widehat S)$, because $\chi(S,\pi^*c\otimes \O_D(D))=0$.  By definition $$\mu(\pi^*L+kD)\otimes E^{\otimes r}=\lambda(\pi^*\beta+k\O_D(D)).$$
Part (2) of Conjecture \ref{blowver} follows directly from the case $k\in [1,\rho-1]$ of  \eqref{psiblow2}. 


\subsection{Blowup formulas for the Segre invariants}
Now we state the blowup relations for the universal functions of the Segre invariants. 

\begin{conjecture}\label{blowsegre1} For $s\in [0,2\rho]$ we have the following formulas.
\begin{enumerate}
\item For $a\in [-\rho,0]$, and furthermore, if $s\in [1,2\rho-1]$, also for  $a\in [-s,-s+\rho]$ we have
\begin{equation}
\label{Sblow1}\sum_{J\subset [\rho-1]} Y_{J,s}^aZ_{J,r}^{-1}= \frac{(1+(2-\smfr{s}{\rho}t^2))^{\frac{1}{2}a(a+s)}}{(1+(1-\smfr{s}{\rho}t^2))^{\frac{1}{2} a(a+s-\rho)}}.\end{equation}
\item For $a\le 0$ we have  
\begin{equation}\label{Sblow2}
\sum_{J\subset [\rho-1]} Y_{J,s}^aZ_{J,s}^{-1}e^{xS_{J,s}}=\frac{(1+(2-\smfr{s}{\rho}t^2))^{\frac{1}{2}a(a+s)}}{(1+(1-\smfr{s}{\rho}t^2))^{\frac{1}{2} a(a+s-\rho)}}e^{\frac{x^2}{2}\cdot t^2(1+(1-\frac{s}{\rho})t^2)+x\cdot at^2 }+O(x^{\rho+a+1}),\end{equation}
furthermore, if $s\in [1,2\rho-2]$ we have for $a\ge -s$ 
\begin{equation}\label{Sblow3}\sum_{J\subset [\rho-1]} Y_{J,s}^aZ_{J,s}^{-1}e^{xS_{J,s}}=\rho \cdot \frac{(1+(2-\smfr{s}{\rho}t^2))^{\frac{1}{2}a(a+s)}}{(1+(1-\smfr{s}{\rho}t^2))^{\frac{1}{2} a(a+s-\rho)}}e^{\frac{x^2}{2}\cdot t^2(1+(1-\frac{s}{\rho})t^2)+x\cdot at^2 }+O(x^{\rho-a-s+1}).\end{equation}
\item Let $\ell\in [1,\rho-1], $ let $a\in \Z$. Then 
\begin{equation}\label{Sblow4}\sum_{J\subset [\rho-1]} \epsilon_\rho^{\ell\|J\|}S^n_{J,s}Y_{J,s}^{a+(\ell-\rho)\frac{s}{\rho}} Z_{J,s}^{-1}=0,\qquad 0\le n<\min\big(\ell(\rho-\ell-a),(\ell+a)(\rho-\ell)\big).\end{equation}
In particular we have
\begin{equation}\label{Sblow5}\sum_{J\subset [\rho-1]} \epsilon_\rho^{\ell\|J\|}Y_{J,s}^{a+(\ell-\rho)\frac{s}{\rho}} Z_{J,s}^{-1}=0,\qquad  a\in [-\ell+1, -\ell+\rho-1].\end{equation}
\end{enumerate}
\end{conjecture}
Like for the Verlinde formulas we want to  see that Conjectures \ref{blowsegre1} and \ref{blowphi} are equivalent, if we assume Conjecture \ref{conj2}. 
We can assume by Conjecture \ref{conj2} that $S$ is a K3 surface, and thus that
 $\phi_{S,\rho,c_1,\alpha,L}\ne 0$.

We put $a=-k$, then, using the definitions of $R_s$ and $Q_s$,  we get by Proposition \ref{verblowprop} the formula
\begin{equation}
\label{phbl}
\begin{split}
\frac{\phi_{\widehat{S},\rho,\pi^*c_1-\ell D,\alpha-a\O_D(D),L-xD}}{\phi_{S,\rho,c_1,\alpha,L}}&=\frac{(1+(1-\smfr{s}{\rho}t^2))^{\frac{1}{2} a(a+s-\rho)}}{(1+(2-\smfr{s}{\rho}t^2))^{\frac{1}{2}a(a+s)}}\\&\quad \cdot e^{-\frac{x^2}{2} \cdot t^2(1+(1-\frac{s}{\rho})t^2)-ax\cdot t^2}\sum_{J \subset [\rho-1]}  \eps_{\rho}^{\|J\| \ell} \frac{Y_{J,s}^{a}}{{Z_{J,s}}}e^{xS_{J,s}}.
\end{split}
\end{equation}
To see that the parts (1) of Conjectures \ref{blowsegre1} and \ref{blowphi}  are equivalent, we put $\ell=0$, $x=0$ in \eqref{phbl}. 

By the case $\ell=0$ of \eqref{phbl},  we see that for any $n>0$ we have $\phi_{\widehat{S},\rho,\pi^*c_1,\alpha-a\O_D(D),L-xD}=\phi_{S,\rho,c_1,\alpha,L}+O(x^n)$, if and only if 
$$\sum_{J\subset [\rho-1]} Y_{J,s}^aZ_{J,s}^{-1}e^{xS_{J,s}}=\frac{(1+(2-\smfr{s}{\rho}t^2))^{\frac{1}{2}a(a+s)}}{(1+(1-\smfr{s}{\rho}t^2))^{\frac{1}{2} a(a+s-\rho)}}e^{\frac{x^2}{2}\cdot t^2(1+(1-\frac{s}{\rho})t^2)+x\cdot at^2 }+O(x^{n}).$$ This shows that
 parts (2) of Conjectures \ref{blowsegre1} and \ref{blowphi} are equivalent.
 
 For part (3), by \eqref{phbl} we just need to see that for a given $n$
we have 
$$\Coeff_{x^m}\Big[e^{-\frac{x^2}{2} \cdot t^2(1+(1-\frac{s}{\rho})t^2)-ax\cdot t^2}\sum_{J \subset [\rho-1]}  \eps_{\rho}^{\|J\| \ell} \frac{Y_{J,s}^{a}}{{Z_{J,s}}}e^{xS_{J,s}}]=0$$ for all $m\le n$, if and only if
$\sum_{J \subset [\rho-1]}S_{J,s}^m  \eps_{\rho}^{\|J\| \ell} \frac{Y_{J,s}^{a}}{{Z_{J,s}}}=0$ for all $m\le n$. This is obvious.

\begin{remark}\label{remeqiv} We can see that  Conjectures \ref{blowver1} and \ref{blowsegre1}  are compatible with the Segre-Verlinde correspondence Conjecture \ref{SegreVerlindeconj}, and the virtual Serre duality Conjecture \ref{SerreDuality}.
\begin{enumerate}
\item
Direct computation gives that replacing $r$ by $-r$ and applying Conjecture \ref{SerreDuality} transforms \eqref{Ablow1} to itself and 
 \eqref{Ablow2} to  \eqref{Ablow2} with $a$ replaced by $-\rho-a$.
 \item 
Direct computation gives that the Segre-Verlinde correspondence Conjecture \ref{SegreVerlindeconj}
transforms \eqref{Sblow1} for $s$ into \eqref{Ablow1} for $r=s-\rho$ and vice versa.
Similarly for $a\in[-\rho+\ell-1,\ell-1]$ it transforms \eqref{Sblow5} for $(a,s)$ into \eqref{Ablow2} for $(a,r=s-\rho)$ and vice versa.
 \end{enumerate}
 \end{remark}

\begin{remark}\label{Bempty}
For every subset $J\subset [\rho-1]$ and all $r,s\in \Z$, using \eqref{Jseries},  we put 
\begin{align*}
\beta_{J,r}&:=\frac{B_{r}}{B_{J,r}}=\frac{1}{\prod_{i\le j\in J} B_{ij,r}},\\
\zeta_{J,s}&:=\frac{Z_{s}}{Z_{J,s}}=\frac{1}{\prod_{i\le j\in J} Z_{ij,s}}.
\end{align*}
Then by definition  $B_{J,r}=\frac{B_r}{\beta_{J,r}}$ and   $Z_{J,r}=\frac{Z_r}{\zeta_{J,r}}$, and the cases $a=0$ of \eqref{Ablow1} and \eqref{Sblow1} give
\begin{align*}
B_{r}&=
\sum_{J\subset [\rho-1]}\beta_{J,r}, \qquad r\in [-\rho,\rho],\\
Z_{s}&=\sum_{J\subset [\rho-1]}\zeta_{J,s}, \qquad s\in [0,2\rho].
\end{align*}
This determines for $r\in [-\rho,\rho]$ the universal function  $B_r$ in terms of the $B_{ij,r}$ and for  $s\in [0,2\rho]$ the universal function $Z_{s}$ in terms of the  $Z_{ij,s}$
\end{remark}

\subsection{Evidence for the conjectures}

In \cite{GKVer} we used Mochizuki's formula \cite[Thm.~7.5.2]{Moc} and localization to compute the Segre and Verlinde invariants of moduli spaces
$M_S^H(\rho,c_1,c_2)$ for 
surfaces $S$ with $p_g(S)>0$ and $b_1(S)=0$, for $\rho=1,2,3,4$, in a number of cases of $S,c_1,c_2$. Assuming 
Conjecture \ref{conj2}, this determines the universal power series $A_{J,r}^{(\rho)}$, $B_{J,r}^{(\rho)}$, $Y_{J,s}^{(\rho)}$, $Z_{J,s}^{(\rho)}$, $S_{J,s}^{(\rho)}$ modulo certain powers of $w$ and $z$. 
Specifically we determined them up to the following orders. 
\begin{enumerate}
\item $A_{J,r}^{(2)}$, $B_{J,r}^{(2)}$, with $r$ as variable modulo $w^{24}$ and $Y_{J,s}^{(2)}$, $Z_{J,s}^{(2)}$, $S_{J,s}^{(2)}$ with $s$ as variable modulo $z^{18}$.
\item $A_{J,r}^{(3)}$, $B_{J,r}^{(3)}$ for $-11\le r\le 3$ modulo $w^{12}$ and $Y_{J,s}^{(3)}$, $Z_{J,s}^{(3)}$, $S_{J,s}^{(3)}$ for $-3\le s \le 12$ modulo $z^{13}$.
\item  Finally we determined $Y_{J,s}^{(4)}$, $Z_{J,s}^{(4)}$, $S_{J,s}^{(4)}$ for $0\le s \le 8$ modulo $z^{8}$.
Using the Verlinde-Segre-correspondence Conjecture \ref{SegreVerlindeconj}, we also determine $A_{J,r}^{(4)}$, $B_{J,r}^{(4)}$ for 
$|r|\le 4$ modulo $w^8$.
\end{enumerate}
Until these orders  Conjectures \ref{blowver1} and \ref{blowsegre1} are confirmed.
In \cite{GKVer} we also conjecturally determined the universal functions  $A_{J,r}^{(\rho)}$, $B_{J,r}^{(\rho)}$, and $Y_{J,s}^{(\rho)}$, $Z_{J,s}^{(\rho)}$, $S_{J,s}^{(\rho)}$, for $\rho=2$, $|r|\le 3$, (respectively $s=-1,\ldots,5$), for $\rho=3$, $|r|\le 3$  (respectively $s=0,\ldots,6$),
and $\rho=4$, $r=-4$ (respectively  $s=0,4$) as algebraic functions. 
We have confirmed modulo $w^{60}$ and $z^{60}$ that theses algebraic functions fulfill Conjectures \ref{blowver1} and \ref{blowsegre1}.

In the current paper we apply some part of these blowup formulas to determine for arbitrary ranks $\rho$ conjectural formulas for the universal functions  $A_{J,r}^{(\rho)}$, $B_{J,r}^{(\rho)}$, and $Y_{J,s}^{(\rho)}$, $Z_{J,s}^{(\rho)}$, $S_{J,s}^{(\rho)}$
for $r=-\rho,0,\rho$ and $s=0$, $\rho$, $2\rho$. We check for $\rho\le 12$ and modulo $w^{60}$ respectively $z^{60}$ that Conjectures \ref{blowver1} and \ref{blowsegre1} hold for these universal functions. In the computations the coefficients are highly overdetermined by the relations. That there is a solution for these relations, which also agrees with the previous results for 
$\rho=1,2,3,4$ is appears quite nontrivial. 

Further confirmation of Conjectures \ref{blowver1} and \ref{blowsegre1} comes from the results of the forthcoming paper \cite{Gstr}, which we describe now.
By Conjecture \ref{poly}  and interpolation, the conjectural formulas for $r=-\rho,0,\rho$ and $s=0,\rho,2\rho$ determine the power series
$A_{J,r}^{(\rho)}$, $B_{J,r}^{(\rho)}$, and $Y_{J,s}^{(\rho)}$, $Z_{J,s}^{(\rho)}$, $S_{J,s}^{(\rho)}$ modulo $w^{3}$ (respectively $z^3$).
Given this information, we use the blowup formulas \eqref{Ablow1}, \eqref{Ablow2}, \eqref{Sblow1}, \eqref{Sblow5}, and the coefficient of $x^1$ of \eqref{Sblow2} and \eqref{Sblow3}
to  compute for $\rho=2,3,4,5,6$ and $|r|\le \rho$ (or $0\le s\le 2\rho$) the coefficients of the power series $A_{J,r}^{(\rho)}$, $B_{J,r}^{(\rho)}$, and $Y_{J,s}^{(\rho)}$, $Z_{J,s}^{(\rho)}$, $S_{J,s}^{(\rho)}$ modulo $w^{100}$ and $z^{100}$ by a Pari/GP program. 
Indeed applying the blowup formulas order by order in $v$ and $t$  gives linear equations for the coefficients of the  $A_{J,r}^{(\rho)}$, $B_{J,r}^{(\rho)}$, and $Y_{J,s}^{(\rho)}$, $Z_{J,s}^{(\rho)}$, $S_{J,s}^{(\rho)}$, which  at each step have a unique solution (indeed they are overdetermined).
The resulting power series fulfill all the relations of Conjectures \ref{blowver1} and \ref{blowsegre1}. Furthermore for $\rho=2,3,4$, they agree up to the order computed with the power series determined in \cite{GKVer} and the algebraic functions conjectured there.

Finally in \cite{Gstr} we use a virtual version of strange duality that relates the power series $A_{J,r}^{(\rho)}$, $B_{J,r}^{(\rho)}$, to the power series  
$A_{J,\rho}^{(r)}$, $B_{J,\rho}^{(r)}$. This allows us to determine from the power series $A_{J,r}^{(\rho)}$, $B_{J,r}^{(\rho)}$ with $|r|\le\rho\le 6$
the $A_{J,r}^{(\rho)}$, $B_{J,r}^{(\rho)}$ with $\rho\le |r|\le 6$. Again we find that the results agree with the results from \cite{GKVer}, both for the 
power series and the algebraic functions conjectured there. In addition with the help of Don Zagier I wrote a program that uses localization on the Hilbert scheme of points to compute
the power series $A_{r}^{(1)}$, $B_{r}^{(1)}$ with $r$ as variable modulo $w^{100}$. The results agree with the power series we 
obtain using the blowup formulas.

The fact that it is possible to determine these power series from the blowup formulas in a consistent way, and whenever  they have been also determined or conjectured by other methods, the results are the same, gives us addition evidence for the validity of the blowup formulas. 

\section{Verlinde formulas for $r=\pm\rho$ and Donaldson invariants in arbitrary rank}
We now will use the blowup formulas  to conjecturally determine  the Verlinde formula for 
$\chi^\vir(M_S^H(\rho,c_1,c_2),\mu(L)\otimes E^{\otimes r})$ in the case $r=-\rho$ and Segre invariants in case $s=0$ for arbitrary rank $\rho$. 
We determine the Donaldson invariants in arbitrary rank $\rho$ as the special case $\alpha=0$ of the Segre invariants for $s=0$.

\subsection{Verlinde formula for $r=-\rho$}

We fix $\rho\in\Z_{>0}$. In order to state the conjectural Verlinde formula for $r=-\rho$,  we introduce the following notations, which we will also use 
in the rest of this paper in the formulation of our main conjectures on Segre and Verlinde invariants.
\begin{notation} Let $\xi=e^{\frac{\pi i}{2\rho}}$ be a primitive $4\rho$-th root of unity.  For $1\le i<j\le \rho-1$ let 
\begin{align*}
\beta_{ij}&:=\frac{\xi^{i+j}-\xi^{-(i+j)}}{\xi^{i-j}-\xi^{j-i}},\quad \beta_{ji}:=\beta_{ij}.
\end{align*}
We put 
\begin{equation} \label{Biij}
\begin{split}
B_{ij}&:=\beta_{ij}^2=\frac{(\xi^{i+j}-\xi^{-i-j})^2}{(\xi^{i-j}-\xi^{-i+j})^2}, \qquad 1\le i<j\le \rho-1,\\
B_{ii}&:=\prod_{j\ne i\in[\rho-1]}\frac{1}{\beta_{ij}}=(-1)^{i-1}\prod_{j\ne i} \frac{\xi^{i-j}-\xi^{-i+j}}{\xi^{i+j}-\xi^{-i-j}},\qquad 1\le i\le \rho-1.
\end{split}
\end{equation}
Then we set
\begin{equation}\label{BJ}
\begin{split}
\beta_I&:=\frac{1}{\prod_{i\le j\in I} B_{ij}},\qquad I\subset [\rho-1],\\
B&:=\sum_{J\subset [\rho-1]} \beta_{J},\\
B_{I}&:=\frac{B}{\beta_{I}}=\sum_{J\subset [\rho-1]} \frac{\beta_{J}}{\beta_I},\qquad I\subset [\rho-1].
\end{split}
\end{equation}
Note that for $I\subset [\rho-1]$ we have 
\begin{align*} \beta_I=\frac{\prod_{i \in I}\prod_{j\ne i}\beta_{ij}}{\prod_{i<j\in I}\beta_{ij}^2}=
\frac{\prod_{i \in I}\prod_{j\ne i}\beta_{ij}}{\prod_{i\ne j\in I}\beta_{ij}}=\prod_{\substack{i\in I\\j\in [\rho-1]\setminus I}}\beta_{ij}.\end{align*}
\end{notation}

Now the Verlinde numbers with $r=-\rho$ are conjecturally determined by the following formulas. 
\begin{conjecture}\label{verrho} Conjecture \ref{conj2} holds for $r=-\rho$ and  any $\rho>0$, if we make  the following definitions.
\begin{enumerate}
\item
$B_{J,-\rho}:=B_J$ for all $J\subset [\rho-1]$, or equivalently $B_{-\rho}:=B$ and 
$B_{ij,-\rho}:=B_{ij}$ for all $1\le i\le j\le \rho-1$.
\item 
\begin{equation}\label{Arel}
A_{-\rho}:=\prod_{j=1}^{\rho-1}\frac{1}{1+\xi^{\rho+2j}w},\quad A_{j,-\rho}:=\frac{1+\xi^{\rho+2j}w}{1+\xi^{\rho-2j}w}, \hbox{ for } j\in[\rho-1],\end{equation}
and therefore for each $J\subset [\rho-1]$
$$
A_{J,-\rho} =\prod_{j\in J}\frac{1}{1+\xi^{\rho-2j}w}\prod_{j\in [\rho-1]\setminus J}\frac{1}{1+\xi^{\rho+2j}w}.
$$
\end{enumerate}
\end{conjecture}

We outline the assumptions, arguments and  calculations that lead to this conjecture. 
We use two working hypotheses.
\begin{enumerate}
\item We expect that for  $J\subset [\rho-1]$ the universal functions $B_{J,-\rho}^{(\rho)}$ are constant, i.e. independent of the variable $w$. 
\item We expect  the formula \eqref{Arel} for the universal power series $A_{J,-\rho}$ to be true.
\end{enumerate}
The motivation for these working hypotheses is the following.
\begin{enumerate}
\item 
Note that, by the Segre-Verlinde correspondence \ref{SegreVerlindeconj}, if $Z_{J,0}^{(\rho)}$ is constant, then  $B_{J,-\rho}^{(\rho)}=Z_{J,0}^{(\rho)}$.
That $Z_{J,0}^{(\rho)}$ is constant is conjectured in \cite{GKVer} as part of a  conjecture about the structure of  the Donaldson invariants in arbitrary rank, and, if one assumes Conjecture \ref{conj2}, it also follows from the Marino-Moore conjecture \cite{MM}. The conjecture of \cite{GKVer} is true for $\rho=1$. For $\rho=2$ and $\alpha=0$ it is the Witten conjecture for algebraic surfaces with $p_g>0$ and $b_1=0$, which was proven in \cite{GNY3}. For $\rho=3$ and $\rho=4$ it was confirmed up to high expected dimension in \cite{GKVer}.

Heuristically we can motivate  that the $Z_{J,0}^{(\rho)}$ should be constant  as follows.
We see from the definition that 
$X_{0}^{(\rho)}=1$ for all $J\subset [\rho-1]$. Thus, if $S$ is a surface with  $p_g(S)>0$ and $b_1(S)=1$ and connected canonical divisor, we get for  $c_1\in H^2(S,\Z)$, that
$$\phi_{S,\rho,c_1,0,0}=\rho^{2-\chi(\O_S)}\sum_{J\subset [\rho-1]}(-1)^{|J|}\eps_{\rho}^{\|J\|K_Sc_1}Z_{J,0}^{K_S^2}.$$
So for $c_2\in \Z$ we get by Conjecture \ref{conj1} that 
$$\Coeff_{z^{\vd(\rho,c_1,c_2)}}\Big[\rho^{2-\chi(\O_S)}\sum_{J\subset [\rho-1]}(-1)^{|J|\chi(\O_S)}\eps_{\rho}^{\|J\|K_Sc_1}Z_{J,0}^{K_S^2}\Big]=\int_{[M_S^H(\rho,c_1,c_2)]^\vir}1,$$
and the right hand side is zero unless $\vd(\rho,c_1,c_2)=0$. This is an indication that  all the $Z_{J,0}$ should be constant.

\item 
For ranks $\rho=2,3,4$ the conjectural formulas for the $A_{J,-\rho}$ were found in \cite{GKVer} based on computations
for not too large virtual dimension, using Mochizuki's formula and localization (for $\rho=4$ also using virtual Serre duality). One can see that these formulas are specializations of  \eqref{Arel}. As the formula is very simple and natural, we expect that it holds for arbitrary rank $\rho$.
\end{enumerate}

For $\rho=2,\ldots 8$, we now define $A_{-\rho}$ and the $A_{j,-\rho}$ by the formulas \eqref{Arel}. We choose variables for the $B_{J,-\rho}$ for $J\subset[\rho-1]$. Imposing the  blowup relations \eqref{Ablow1}, \eqref{Ablow2}  order by order in $v=w$ gives an infinite number of linear equations for the $B_{J,-\rho}$, which we can solve successively using Pari/GP. We find that the $B_{J,-\rho}$ are uniquely determined by these equations, (and that they fulfill them to high order in $w$). We analize the explicit formulas that we obtain for $B_{-\rho}$ and the $B_{ij,-\rho}$ for  $\rho=2,\ldots 8$. Using also Remark \ref{Bempty}, we find that, for $\rho=2,\ldots,8$, we can write 
$$\frac{B_{\emptyset,-\rho}}{B_{J,-\rho}}=\prod_{\substack{i\in I\\j\in [\rho-1]\setminus I}}\beta_{ij},$$
 for suitable numbers $\beta_{ij}=\beta_{ji}$, which we finally see as given by formula \eqref{Biij}. This gives the formulas of Conjecture \ref{verrho} for $\rho=2,\ldots,8$.
Using Pari/GP we check that  these  universal functions satisfy  the blowup formulas \eqref{Ablow1} and \eqref{Ablow2}  as identities of rational functions in $w$ for $\rho\le 7$ and $\rho\le 6$  respectively, and that they satisfy   \eqref{Ablow1}, \eqref{Ablow2} modulo $O(w^{100})$  for $\rho\le 12$.

If we chose $\rho,c_1,c_2$ with $\vd(\rho,c_1,c_2)=0$, then the degree of the virtual fundamental class $\int_{[M_S^H(\rho,c_1,c_2)]^\vir} 1$, is both 
the coefficient of  $w^0$ of $\psi_{S,\rho,c_1,L,r}$ for all $r\in \Z$ (using the virtual Riemann-Roch formula \cite{FG})  and the coefficient of  $z^0$ of $\phi_{S,\rho,c_1,\alpha,L}$ 
for any $\alpha\in K^0(S)$. 
Thus  we conjecture that the constant coefficient of each $Z_{J,s}$ and $B_{J,r}$ for all $r,s\in \Z$ is $B_J$. Using also Conjecture \ref{poly} we get the following statement.
\begin{conjecture}
\begin{enumerate}
\item For all $J\subset [\rho-1]$ and $1\le i\le j\le \rho-1$ and all $r\in \Z$, we  can write $B_{J,r}=B^0_{J,r}B_J$, $B_{r}=B^0_rB$, $B_{ij,r}=B^0_{ij,r}B_{ij}$ with $B^0_{J,r}$, $B^0_r$, $B^0_{ij,r}\in 1+w\C[[w]]$.
\item In the same way for all
$J\subset [\rho-1]$ and  $1\le i\le j\le \rho-1$ and all $r\in \Z$, we we can write $Z_{J,s}=Z^0_{J,s}B_J$, $Z_{r}=Z^0_rB$, $Z^0_{ij,r}=Z^0_{ij,r}B_{ij}$ with $Z^0_{J,r}$, $Z^0_r$, $Z^0_{ij,r}\in 1+z\C[[z]$.
\end{enumerate}
\end{conjecture}

\subsection{Segre numbers for $s=0$ and Donaldson invariants in arbitrary rank}
Applying the Segre-Verlinde correspondence Conjecture \ref{SegreVerlindeconj} to Conjecture \ref{verrho} determines the universal functions $Y_{J,0}$, $Z_{J,0}$. Thus in order to get the Segre numbers for $s=0$ in arbitrary rank we only need to determine the $S_{J,0}$. We obtain the following conjecture, 
which generalizes the results of \cite{GKVer} for $\rho\le 4$.
\begin{conjecture}\label{Segre0} Under the change of variables $z=t(1+t^2)^{\frac{1}{2}}$, Conjectures \ref{conj1} \and \ref{conj2} hold with 
\begin{align*}
Z_{J,0}&:=B_{J,0}\quad \hbox{ for all $J\subset [\rho-1]$,}\\
Y_{J,0}&:=\frac{(1+t^2)^\rho}{(1+2t^2)^{\frac{1}{2}}}\prod_{j\in J} \frac{1}{(1+t^2)^{\frac{1}{2}}+\xi^{\rho-2j}t}\prod_{j\in [\rho-1]\setminus J} \frac{1}{(1+t^2)^{\frac{1}{2}}+\xi^{\rho+2j}t},\\
S_{J,0}&:=\Big(-\sum_{j\in J}\xi^{\rho-2j}+\sum_{j\in [\rho-1]\setminus J}\xi^{2j-\rho}\Big)z.
\end{align*}
In particular, we have 
\begin{equation}\label{S0form}
S_{0}=\Big(\sum_{j\in [\rho-1]}\xi^{2j-\rho}\Big)z,\quad  S_{j,0}=(-\xi^{2j-\rho}-\xi^{\rho-2j})z=-2\sin(\pi\frac{j}{\rho})z.\end{equation}
\end{conjecture}
To determine the $S_{J,0}$, we make the working assumption that $S_{0}$ and the  $S_{i,0}$  for $i=1,\ldots, \rho-1$ are constant multiples  $s_{0} z$, $s_{i,0}z$ of $z$ with $s_0\in \C$, $s_{i,0}\in \C$. For $\rho=2,\ldots,8$ we first let $s_0$ and the $s_{i,0}$ be variables. Then the coefficients of $x^1$ of the  blowup formulas \eqref{Sblow2} impose, order by order in $t$, linear equations for $s_0$ and the $s_{i,0}$, which uniquely determine them by the formula of \eqref{S0form}. 
Then for $\rho\le 12$, we define $S_{0}$ and the  $S_{i,0}$ by \eqref{S0form} and check using Pari/GP, that for $\rho\le 12$ the blowup formulas \eqref{Sblow1}, \eqref{Sblow2},  \eqref{Sblow3}, \eqref{Sblow4} hold modulo $O(t^{100})$.

Conjecture \ref{Donconj} follows from Conjecture \ref{Segre0} by specializing to $\alpha=0$:
By Definition $D^{S,H}_{\rho,c_1,c_2}(L  +u \pt)=\Coeff_{z^\vd(\rho,c_1,c_2)}\big[\phi_{S,\rho,c_1,0,L}\big]$.
As in the introduction we put $\widetilde a:=2a-K_S$ for $a\in H^2(S,\Z)$, and write $\widetilde{SW}(\widetilde a):=SW(a)$.
Note that $$\frac{1}{2}\sum_{j=1}^{\rho-1}(\xi^{2j-\rho}+\xi^{\rho-2j})=\frac{1}{2}\sum_{j=1}^{\rho-1}(\xi^{2j-\rho}+\xi^{2(\rho-j)-\rho})=\sum_{j=1}^{\rho-1}\xi^{2j-\rho},$$
therefore $$\prod_{j=1}^{\rho-1}e^{-\sin(\pi\frac{j}{\rho})(\widetilde a_j L)z}=e^{(K_SL)S_{0}z}\prod_{j=1}^{\rho-1}e^{(a_j L)S_{j,0}z}.$$
Using $a_iK_S=a_i^2$, $a_jK_S=a_j^2$ we see that $\frac{1}{2}\widetilde a_i(\widetilde a_j-\widetilde a_i)=2a_ia_j-a_i^2-a_j^2$.  By $B_{ij}=\beta_{ij}^2$ and 
$B_{ii}=\prod_{j\in [\rho-1]\setminus \{i\}} \beta_{ij}^{-1}$, 
this gives 
$$\prod_{1\le i< j\le \rho-1}\beta_{ij}^{\frac{1}{2}\widetilde a_i(\widetilde a_j-\widetilde a_i)}=\prod_{1\le i< j\le \rho-1} B_{ij}^{a_ia_j}\cdot  \prod_{i\ne j}\beta_{ij}^{-a_i^2}=\prod_{1 \leq i \leq j \leq \rho-1} B_{ij}^{a_i a_j}.$$
Therefore Conjecture \ref{Donconj} follows from Conjecture \ref{Segre0} and Conjecture \ref{conj2}.

\begin{remark}
The Donaldson invariants are more generally defined for differentiable $4$-manifolds $M$ with $b_{+}(M)>1$ and $b_1(M)=0$. If $M$ is a complex projective surface, they agree up to different sign conventions with the definition of the Donaldson invariants above. Let $\sigma(M)=b_+(M)-b_-(M)$ be signature of  the intersection form on $H_2(M,\Q)$, and 
$e(M)$ the topological Euler characteristic of $M$. Replacing in \eqref{Donform} $K_S^2$ by $3\sigma(M)+2e(M)$ and  $\chi(\O_S)$ by $\smfr{1}{4}(\sigma(M)+e(M))$, we obtain a formula that (up to sign) is expressed in terms of invariants of the $4$-manifold $M$. 
Note that $\frac{c_1^2-c_1K_S}{2}\in \Z$. 
An easy calculation shows that  changing the sign from $\prod_{j=1}^{\rho-1}\varepsilon_\rho^{j\cdot(a_j c_1)}$ to $\prod_{j=1}^{\rho-1}\varepsilon_\rho^{j\cdot\frac{1}{2}(\widetilde a_j+c_1)c_1}$, 
just multiplies \eqref{Donform} by the global sign $(-1)^{(\rho-1)\frac{c_1^2-c_1K_S}{2}}$. 

We conjecture that, after making these replacements, formula \eqref{Donform} gives  (up to possibly different sign conventions) the Donaldson invariants of rank $\rho$ of any simply connected $4$ manifold $M$ of Seiberg-Witten simple type.
\end{remark}

\subsection{Verlinde formula for $r=\rho$ and Segre formula for $s=2\rho$}
Using  the virtual Serre duality Conjecture \ref{SerreDuality}, we immediately get  from Conjecture \ref{verrho} a Verlinde formula for the case $r=\rho$.
In addition we checked using Pari/GP that  the universal functions $A_{J,\rho}$, $B_{J,\rho}$ below satisfy the blowup formulas \eqref{Ablow1}, \eqref{Ablow2} as identities of rational functions in $w$ for $\rho\le 6$ and $\rho\le 5$ respectively, and that they satisfy   \eqref{Ablow1}, \eqref{Ablow2} modulo $O(w^{100})$  for $\rho\le 12$.
\begin{conjecture}\label{vermrho} Conjecture \ref{conj2} holds for $r=\rho$  with  for all $J\in [\rho-1]$
\begin{align*}A_{J,\rho}&=(1+w^2)^{1-\rho}\prod_{j\in J}  (1-\xi^{\rho-2j}w)\prod_{j\in [\rho-1]\setminus J}(1-\xi^{\rho+2j}w),\\
B_{J,\rho}&=(1+w^2)^{\binom{\rho}{2}}B_{J}\prod_{j\in J}  \frac{1}{(1-\xi^{\rho-2j}w)^\rho}\prod_{j\in [\rho-1]\setminus J}\frac{1}{(1-\xi^{\rho+2j}w)^\rho}.
\end{align*}
\end{conjecture}
We can apply the Segre-Verlinde correspondence to Conjecture \ref{vermrho}, and then use the blowup formulas again to get the Segre invariants for $s=2\rho$. Again the results agree  with those of \cite{GKVer} for $\rho\le 3$.
\begin{conjecture}\label{Serre2rho}
Conjecture \ref{conj2} holds for $s=2\rho$, when we put for $J\subset[\rho-1]$
\begin{align*}
Y_{J,2\rho}&=(1+z^2)^{\frac{1-\rho}{2}}\prod_{j\in J}(1-\xi^{\rho-2j}z)\prod_{j\in [\rho-1]\setminus J}(1-\xi^{\rho+2j}z),\\
Z_{J,2\rho}&=(1+z^2)^{\binom{\rho}{2}}B_{J}\prod_{j\in J}  \frac{1}{(1-\xi^{\rho-2j}z)^\rho}\prod_{j\in [\rho-1]\setminus J}\frac{1}{(1-\xi^{\rho+2j}z)^\rho},\\
S_{J,2\rho}&=\frac{\rho+1}{1+z^2}-\frac{2\rho}{1-(-z^2)^\rho}+\sum_{j\in J}\frac{1}{1-\xi^{\rho+2j}z}+ \sum_{j\in [\rho-1]\setminus J}\frac{1}{1-\xi^{\rho-2j}z}.\\
\end{align*}
\end{conjecture}
The formulas for the $Y_{J,\rho}$ and $Z_{J,\rho}$ follow  from  Conjecture \ref{vermrho} and the Segre-Verlinde correspondence Conjecture \ref{SegreVerlindeconj}.
In fact, note that we have for $r=\rho$ the variable change $w=v$, and in the Segre-Verlinde correspondence we have the change of variables $v=t(1-t^2)^{-\frac{1}{2}}$, and finally we have 
$z=t(1-t^2)^{-\frac{1}{2}}$, i.e. $t=z/(1+z^2)^{\frac{1}{2}}$. Combining these, we get the trivial change of variables $w=z$, and we have $W_{2\rho}=(1-t^2)^{\frac{1-\rho}{2}}=(1+z^2)^{\frac{1-\rho}{2}}.$

For $\rho=2,\ldots,6$ we let $S_{2\rho}$ and the $S_{i,2\rho}$ for $i=1,\ldots,5$ be power series in $z$ with vanishing constant terms and independent variables as coefficients. We impose order by order in $z$ the coefficient of $x^1$ of the blowup formulas
\eqref{Sblow2}, \eqref{Sblow4}. These impose linear relations on the coefficients of $S_{2\rho}$ and the $S_{i,2\rho}$, which we check have a unique solution 
modulo $O(z^{100})$, given by the above formulas. We then check with a Pari/GP program, that with these formulas  Conjecture 3.3. holds for $\rho\le 12$ modulo $O(z^{100})$.

\section{ $K$-theoretic Donaldson invariants and Segre invariants for $s=\rho$}
In this section we give a conjectural formula for the holomorphic Euler characteristics $\chi^\vir(M_S^H(\rho,c_1,c_2),\mu(L))$ for arbitrary
rank $\rho>0$. 
These are also called the $K$-theoretic Donaldson invariants, and have been studied in case $\rho=2$ e.g. in \cite{GNY2},\cite{GKW},\cite{GY},\cite{Got1}. 
Finally we use the Verlinde-Segre correspondence and the blowup formulas to also determine the  Segre invariants with $s=\rho$. 
 
  \subsection{K-theoretic Donaldson invariants }
The universal functions for the $K$-theoretic Donaldson invariants are determined by the following formulas, which generalize the formulas of \cite{GKVer} for $\rho\le 3$. These formulas  have a number of similarities with the formulas in case for $r=\pm\rho$.

 \begin{conjecture} \label{ver0} Let again $\xi=e^{\frac{\pi i}{2\rho}}$. Then we have the following.
 \begin{align}
 \label{A0}\begin{split}A_{0}&=\prod_{i=1}^{\rho-1}\Big(\big(1+\big(\smfr{\xi^{\rho-2i}-\xi^{2i-\rho}}{2})^2w^2\big)^{\frac{1}{2}}-\big(\smfr{\xi^{\rho-2i}+\xi^{2i-\rho}}{2})w\Big),\\
 A_{i,0}&=\frac{\big(1+\big(\smfr{\xi^{\rho-2i}-\xi^{2i-\rho}}{2})^2w^2\big)^{\frac{1}{2}}+\big(\smfr{\xi^{\rho-2i}+\xi^{2i-\rho}}{2})w}{\big(1+\big(\smfr{\xi^{\rho-2i}-\xi^{2i-\rho}}{2})^2w^2\big)^{\frac{1}{2}}-\big(\smfr{\xi^{\rho-2i}+\xi^{2i-\rho}}{2})w},\qquad i=1,\ldots, \rho-1,
  \end{split}\\
 \label{gamij}
 \begin{split}
 \gamma_{i,\rho-i}&=\Big(\frac{1+\big(\smfr{\xi^{\rho-2i}-\xi^{2i-\rho}}{2})^2w^2}{1-w^2}\Big)^{\frac{1}{2}},\qquad\qquad i=1,\ldots,\rho-1,\\
 \gamma_{i,j}&=\frac{(\xi^{\rho-2j}+\xi^{2j-\rho})\gamma_{i,\rho-i}+(\xi^{\rho-2i}+\xi^{2i-\rho})\gamma_{j,\rho-j}}{\xi^{\rho-2i}+\xi^{2i-\rho}+\xi^{\rho-2j}+\xi^{2j-\rho}},\quad 1\le i\ne j\le \rho-1,
 \end{split}\\
\label{gamI} \gamma_{I}&=\Big(\prod_{i\in I,\ j\not\in I}  \gamma_{ij}\Big) \Big(\prod_{i\in I} A_{i,0}\Big)^{\frac{\rho}{2}}, \qquad \qquad I\subset [\rho-1],\\
 \label{BJ0}
 B_0&=\sum_{I\subset [\rho-1]} \beta_I\gamma_I,\qquad B_{J,0}=\frac{B_0}{\beta_J\gamma_J},\qquad J\subset[\rho-1].
  \end{align}
 \end{conjecture}

 We explain the strategy for finding these formulas: we use the virtual Serre duality relation. In geneal it relates $r$ to $-r$, but in our case  $r=0$ it gives nontrivial relations among the  generating functions $A_{J,0}$, $B_{J,0}$. It implies for $J\in [\rho-1]$ the relations
\begin{align}\label{0A} A_{J,0}(w)&=\frac{(1-w^2)^{\rho-1}}{A_{J,0}(-w)},\\
\label{0B} B_{J,0}(w)&=\frac{A_{J,0}(-w)^\rho B_{J,0}(-w)}{(1-w^2)^{\binom{\rho}{2}}}.
\end{align}
In particular we have the weaker relation
\begin{equation} \label{0B1} \frac{B_{J,0}(w)}{B_0(w)}=\frac{B_{J,0}(-w)}{B_0(-w)}\prod_{i\in J}A_{i,0}(-w)^\rho .
\end{equation}
In \cite{GKVer} the universal series $A_{J,0}$, $B_{J,0}$ have been determined for $\rho\le 3$. 
By Conjectures \ref{verrho}, \ref{vermrho}, we have conjectural formulas for  the $A_{J,r}$, $B_{J,r}$ when $r=\pm \rho$.
Using Conjecture \ref{poly}, this allows us to use interpolation to  determine the $A_{J,r}$, $B_{J,r}$ up to degree $1$ in $w$. Now for 
$\rho\le 4$ we take the coefficients of the power series $A_{0}$, $A_{i,0}$, $B_{0}$, $B_{ij,0}$ as indeterminants. Then the blowup relations 
\eqref{Ablow1}, \eqref{Ablow2} give successively degree by degree in $w$ linear relations for these coefficients. We use them to determine the power series 
modulo $O(w^{60})$ with a Pari/GP program. From these truncated power series and \eqref{0A} we guess the formulas \eqref{A0}  for $\rho\le 4$. We then make the Ansatz that  the formula \eqref{A0} extends to arbitrary rank $\rho$. Note that this makes \eqref{0A} automatic.

In the case $r=-\rho$ we found for $J\subset[\rho-1]$ that 
$$B_{-\rho}=\sum_{I\subset [\rho-1]} \beta_I,\qquad B_{J,-\rho}=\frac{B_{-\rho}}{\beta_J}, \qquad \beta_{J}=\prod_{i\in J\ j\not\in J}\beta_{ij},$$
where $\beta_{ij}=\beta_{ji}$ are constants.  
We expect the formula for the $B_{J,0}$ to have a similar structure. We make the Ansatz that 
\eqref{BJ0} and \eqref{gamI} hold
 with unknown power series $\gamma_{ij}=\gamma_{ji}\in 1+w^2\C[[w^2]]$. This  is also motivated both by Remark \ref{Bempty}, and the fact that
 with this Ansatz at least the relation \eqref{0B1} is automatically fullfilled. In fact this appears to be the simplest modification of formula in case $r=-\rho$ compatible with  Remark \ref{Bempty} and  \eqref{0B1}. The relation \eqref{0B} will not automatically follow from the Ansatz but requires the right choice of the $\gamma_{ij}$.
 
 Now for $\rho\le 8$, we take the coefficients of positive degree of the $\gamma_{ij}$ for $1\le i<j\le \rho-1$ as indeterminates. Then the blowup relations \eqref{Ablow1}, \eqref{Ablow2} give, successively degree by degree in $w$, linear relations for these coefficients, which we use with a Pari/GP program to determine them modulo $O(w^{40})$. We analize these power series, and find  the formula for the $\gamma_{i,\rho-i}$, and then that the general $\gamma_{ij}$ are given by the convex combination \eqref{gamij}.
 We then use a Pari/GP program to show that with these formulas  the formula \eqref{0B} and the blowup formulas Conjecture  \ref{blowver1} hold modulo $O(w^{80})$ for $\rho \le 10$.

 \subsection{Segre invariants for $s=\rho$}
 We can now apply the Segre Verlinde correspondence to determine the universal power series $Y_{J,\rho}$, $Z_{J,\rho}$ of the 
 Segre invariants for $s=\rho$. Then the power series $S_{J,\rho}$ can be determined using the blowup formulas of Conjecture \ref{blowsegre1}.
 We obtain the following  formulas.

  \begin{conjecture} \label{segrerho} Let again $\xi=e^{\frac{\pi i}{2\rho}}$. We put 
 \begin{align*}
  Y_{\rho}&=\prod_{i=1}^{\rho-1}\Big(\big(1+\big(\smfr{\xi^{\rho-2i}+\xi^{2i-\rho}}{2})^2z^2\big)^{\frac{1}{2}}+\big(\smfr{\xi^{\rho-2i}+\xi^{2i-\rho}}{2})z\Big),\\
Y_{i,\rho}&= \Big(\big(1+\big(\smfr{\xi^{\rho-2i}+\xi^{2i-\rho}}{2})^2z^2\big)^{\frac{1}{2}}-\big(\smfr{\xi^{\rho-2i}+\xi^{2i-\rho}}{2})z\Big)^2\\ \zeta_{i,\rho-i}&=\big(1+\big(\smfr{\xi^{\rho-2i}+\xi^{2i-\rho}}{2})^2z^2\big)^{\frac{1}{2}},\qquad\qquad i=1,\ldots,\rho-1,\\
 \zeta_{i,j}&=\frac{(\xi^{\rho-2j}+\xi^{2j-\rho})\zeta_{i,\rho-i}+(\xi^{\rho-2i}+\xi^{2i-\rho})\zeta_{j,\rho-j}}{\xi^{\rho-2i}+\xi^{2i-\rho}+\xi^{\rho-2j}+\xi^{2j-\rho}},\quad 1\le i\ne  j\le \rho-1,\\
 \zeta_{I}&=\Big(\prod_{i\in I,\ j\not\in I}  \zeta_{i,j}\Big)\Big(\prod_{i\in I} Y_{i,\rho}^{\frac{\rho}{2}}\Big),\qquad \qquad I\subset [\rho-1],\\
 Z_{\rho}&=\sum_{I\subset [\rho-1]} \beta_I\zeta_I, \qquad Z_{J,\rho}=\frac{Z_{\rho}}{\beta_J\zeta_J},\qquad J\subset[\rho-1],\\
   S_{\rho}&=-\sum_{i=1}^{\rho-1}\Big(\big(\smfr{\xi^{\rho-2i}+\xi^{2i-\rho}}{2}\big)^2z^2+\big(\smfr{\xi^{\rho-2i}+\xi^{2i-\rho}}{2})z\big(1+\big(\smfr{\xi^{\rho-2i}+\xi^{2i-\rho}}{2})^2z^2\big)^{\frac{1}{2}}\Big),\\
 S_{i,\rho}&=-\smfr{\xi^{2(\rho-2i)}-\xi^{2(2i-\rho)}}{2}z^2-\big(\xi^{\rho-2i}+\xi^{2i-\rho}\big)z\big(1+\big(\smfr{\xi^{\rho-2i}+\xi^{2i-\rho}}{2})^2z^2\big)^{\frac{1}{2}},\quad i=1,\ldots,\rho-1.
 \end{align*}
 \end{conjecture}
 
The formulas for the $Y_{J,\rho}$ and the $Z_{J,\rho}$ follow directly from Conjecture \ref{ver0}.
We let $S_\rho$, and $S_{i,\rho}$ for $i=1,\ldots,\rho-1$ be power series in $z$ with constant coefficient $0$ and the other coefficients as variables. Then the coefficients of $x^1$ of the blowup formulas \eqref{Sblow2} give linear relations for the coefficients of the 
$S_{i,\rho}$ which we use to determine them up to degree $50$ in $z$ for $\rho \le 7$. 
From these we can read off the formulas for  $S_{\rho}$ and the $S_{i,\rho}$. We then check that with the formulas of Conjecture \ref{segrerho},  the blowup formulas of Conjecture \ref{blowphi} hold up to degree $60$ in $z$.

{\tt{gottsche@ictp.it}}

\begin{thebibliography}{AGDP}
\bibitem[Don]{Don} S.K.~Donaldson, \textit{Polynomial  invariants  for  smooth  four-manifolds}, Topol.~29 (1990) 257--315.
\bibitem[EGL]{EGL} G.~Ellingsrud, L.~G\"ottsche, and M.~Lehn, \textit{On the cobordism class of the Hilbert scheme of a surface}, Jour.~Alg.~Geom.~10 (2001) 81--100. 
\bibitem[FG]{FG} B.~Fantechi and L.~G\"ottsche, \textit{Riemann-Roch theorems and elliptic genus for virtually smooth schemes}, Geom.~Topol.~14 (2010) 83--115.
\bibitem[FL1]{FL1} P.~Feehan and T.~Leness,  \textit{The $\mathrm{SO}(3)$ monopole cobordism and superconformal simple type}, Adv.~Math.~356 (2019) 106817.
\bibitem[FL2]{FL2} P.~Feehan and T.~Leness,  \textit{Superconformal simple type and Witten's conjecture}, Adv.~Math.~356 (2019) 106821.

\bibitem[G1]{Got1} L.~G\"ottsche, \textit{Verlinde-type formulas for rational surfaces}, Journal of the EMS, {\bf 22} (2019), 151--212.
\bibitem[G2]{Gstr} L.~G\"ottsche, \textit{Blowup formula, strange duality and the Verlinde and Segre numbers of elliptic surfaces}, in preparation.
\bibitem[GKL]{GKL} L.~G\"ottsche, M.~Kool and T.~Laarakker, \textit{Virtual degeneracy loci, Ramanujan's continued fractions, and cosmic strings}, arXiv:2108.13413 .
\bibitem[GKVer]{GKVer} L.~G\"ottsche and M.~Kool, \textit{Virtual Segre and Verlinde numbers of projective surfaces}, arXiv:2007.11631.
\bibitem[GKW]{GKW}  L. G\"ottsche, M.~Kool, and R.A.~Williams, \textit{Verlinde formulae on complex surfaces: K-theoretic invariants}, Forum of Math.~Sigma 9 (2021) 1--31
\bibitem[GNY1]{GNY2} L.~G\"ottsche, H.~Nakajima, and K.~Yoshioka, \textit{K-theoretic Donaldson invariants via instanton counting}, Pure and Appl.~Math.~Quart.~5 (2009) 1029--1111.
\bibitem[GNY2]{GNY3} L.~G\"ottsche, H.~Nakajima, and K.~Yoshioka, \textit{Donaldson = Seiberg-Witten from Mochizuki's formula and instanton counting}, Publ.~Res.~Inst.~Math.~Sci.~47 (2011) 307--359.
\bibitem[GY]{GY} L.~G\"ottsche and Y.~Yuan, \textit{Generating functions for $K$-theoretic Donaldson invariants and Le Potier's strange duality}, J.~Alg.~Geom.~28 (2019) 43--98. 
\bibitem[Joh]{Joh} D.~Johnson, \textit{Universal series for Hilbert schemes and strange duality}, IMRN 2020 10 (2020) 3130--3152.
\bibitem[Kro]{Kro} P.~B.~Kronheimer, \textit{Four-manifold invariants from higher-rank bundles}, Jour.~Diff.~Geom.~70 (2005) 59--112.
\bibitem[Leh]{Leh} M.~Lehn, \textit{Chern classes of tautological sheaves on Hilbert schemes of points on surfaces}, Invent.~Math.~136 (1999) 157--207.
\bibitem[LP]{LP} Le Potier, J.: \textit{ Dualit\'e \'etrange, sur les surfaces}. Preliminary version 18.11.05
\bibitem[MM]{MM} M.~Mari\~{n}o and G.~Moore, \textit{The  Donaldson-Witten function for gauge groups of rank larger than one}, Comm.~Math.~Phys.~199 (1998) 25--69.
\bibitem[MOP1]{MOP1} A.~Marian, D.~Oprea, and R.~Pandharipande, \textit{Segre classes and Hilbert schemes of points}, Ann.~Sci.~ENS 50 (2017) 239--267.
\bibitem[MOP2]{MOP2} A.~Marian, D.~Oprea, and R.~Pandharipande, \textit{The combinatorics of Lehn's conjecture}, J.~Math.~Soc.~Japan 71 (2019) 299--308.
\bibitem[MOP3]{MOP3} A.~Marian, D.~Oprea, and R.~Pandharipande, \textit{Higher rank Segre integrals over the Hilbert scheme of points}, arXiv:1712.02382.
\bibitem[Moc]{Moc} T.~Mochizuki, \textit{Donaldson type invariants for algebraic surfaces}, Lecture Notes in Math.~1972, Springer-Verlag, Berlin (2009). 
\bibitem[Mor]{Mor} J.~W.~Morgan, \textit{The Seiberg-Witten equations and applications to the topology of smooth four-manifolds}, Math.~Notes 44, Princeton Univ.~Press (1996).
\bibitem[Voi]{Voi} C.~Voisin, \textit{Segre classes of tautological bundles on Hilbert schemes of surfaces}, Alg.~Geom.~6 (2019) 186--195.
\bibitem[Wit]{Wit} E.~Witten, \textit{Monopoles and four-manifolds}, Math.~Res.~Lett.~1 (1994) 769--796. 
\end{thebibliography}
\end{document}